\documentclass[12pt]{amsart}
\usepackage{latexsym,amssymb,amsrefs}

\newtheorem{lemma}{Lemma}

\newtheorem{theorem}[lemma]{Theorem}
\newtheorem{corollary}[lemma]{Corollary}
\newtheorem{definition}[lemma]{Definition}

\theoremstyle{definition}

\renewcommand{\MR}{\ }

\newcommand{\proofend}{\begin{flushright}
                       $\Box$
                       \end{flushright}}

\begin{document}

\title[Morse-Bott inequalities]{The Morse-Bott inequalities via dynamical systems}

\author{Augustin Banyaga}
\address{Department of Mathematics \\
         Penn State University \\
         University Park, PA 16802}
\email{banyaga@math.psu.edu}

\author{David E. Hurtubise}
\address{Department of Mathematics and Statistics\\
         Penn State Altoona\\
         Altoona, PA 16601-3760}
\email{Hurtubise@psu.edu}


\begin{abstract}
Let $f:M \rightarrow \mathbb{R}$ be a Morse-Bott function on a compact smooth
finite dimensional manifold $M$. The polynomial Morse inequalities and an explicit 
perturbation of $f$ defined using Morse functions $f_j$ on the critical
submanifolds $C_j$ of $f$ show immediately that $MB_t(f) = P_t(M) + (1+t)R(t)$, where
$MB_t(f)$ is the Morse-Bott polynomial of $f$ and $P_t(M)$ is the Poincar\'e 
polynomial of $M$. We prove that $R(t)$ is a polynomial with nonnegative 
integer coefficients by showing that the number of gradient flow lines
of the perturbation of $f$ between two critical points $p,q \in C_j$ of relative
index one coincides with the number of gradient flow lines between $p$ and $q$ of
the Morse function $f_j$. This leads to a relationship between the kernels of the 
Morse-Smale-Witten boundary operators associated to the Morse functions $f_j$ and 
the perturbation of $f$. This method works when $M$ and all the critical submanifolds
are oriented or when $\mathbb{Z}_2$ coefficients are used.
\end{abstract}

\maketitle


\section{Introduction}

Let $h: M \rightarrow  \mathbb{R}$ be a Morse function on a compact smooth
manifold of dimension $m$. The Morse inequalities say that
$$
\nu_n  - \nu_{n-1} + \cdots + (-1)^n \nu_0 \geq  b_n - b_{n-1} + \cdots +(-1)^n b_0
$$
for all $n = 0,\ldots, m$ (with $\geq$ an equality when $n=m$), where $\nu_k$ is 
the number of critical points of $h$ of index $k$ and $b_k$ is the $k^\text{th}$
Betti number of $M$ for all $k$. These inequalities follow from the fact 
that the Morse function $h$ determines a CW-complex $X$ whose cellular homology
is isomorphic to the singular homology of $M$ (the $k$-cells of $X$ are in bijective
correspondence with the critical points of index $k$).

\vspace{.1 in}

In \cite{WitSup} Witten introduced the idea that the Morse inequalities can be studied by
deforming the de Rham differential on differential forms by the differential of the Morse
function. This led him to consider a chain complex whose chains are generated by the 
critical points of the Morse function and whose differential is defined by counting the gradient
flow lines between critical points of relative index one (the so-called Morse-Smale-Witten chain
complex). See also \cite{BotMor}. The homology of this complex is called the ``Morse homology'',
and the Morse Homology Theorem asserts that the Morse homology is isomorphic to the singular
homology (see \cite{BanLec} and \cite{SchMor}). For an excellent exposition of Witten's ideas see \cite{HenLes}.

On one hand Bismut, and on the other hand, Helffer and Sj\"ostrand have given rigorous
mathematical derivations of Witten's analytical ideas \cite{BisThe}\cite{HelApr}.
In this way, they proved the Morse inequalities and the Morse-Bott inequalities, which relate
the Betti numbers of $M$ and the Betti numbers of the critical submanifolds of a 
Morse-Bott function on $M$, without using the Morse Homology Theorem.

\vspace{.1 in}

In this paper, we present a proof of the Morse-Bott inequalities using ideas from
dynamical systems. The proof uses the Morse inequalities, the Morse Homology Theorem,
and an explicit perturbation technique that produces a Morse-Smale function $h$ arbitrarily 
close to a given Morse-Bott function $f$ \cite{AbbLec} \cite{AusMor}. Roughly, the
perturbation $h$ of $f$ is constructed as follows.  We first fix a Riemannian metric
on $M$ and apply the Kupka-Smale Theorem to choose Morse-Smale functions $f_j$
on the critical submanifolds $C_j$ for $j=1,\ldots ,l$. Next, we extend the Morse-Smale
functions $f_j$ to tubular neighborhoods $T_j$ of the critical submanifolds, and we define
$$
h = f + \varepsilon \left( \sum_{j=1}^l \rho_j f_j \right)
$$
where $\rho_j$ is a bump function on $T_j$ and $\varepsilon > 0$.
We then apply a well-known folk theorem (whose proof can be found in
Section 2.12 of \cite{AbbLec}) which says that we can perturb the Riemannian metric
on $M$ outside of the union of the tubular neighborhoods $T_j$ for $j=1,\ldots ,l$
so that $h$ satisfies the Morse-Smale transversality condition with respect to the
perturbed metric.

Once we know that $h$ and $f_j$ for $j=1,\ldots ,l$ satisfy the Morse-Smale transversality
condition, we can compare the Morse-Smale-Witten chain complex of $h$ to the
Morse-Smale-Witten chain complexes of $f_j$ for $j=1,\ldots ,l$. This allows us
to show that the coefficients of the polynomial $R(t)$ in Theorem \ref{MorseBottineq}
are nonnegative. The proof works when the manifold $M$ and all the critical submanifolds
are oriented or when $\mathbb{Z}_2$ coefficients are used. 

\nocite{BotNon}
\nocite{BotLec}
\nocite{BanMor}
\nocite{JiaMor}


\section{The Morse-Smale-Witten chain complex}\label{MorseChain}

In this section we briefly recall the construction of the Morse-Smale-Witten
chain complex and the Morse Homology Theorem.  For more details see \cite{BanLec}.

\smallskip
Let $Cr(f) = \{p \in M |\, df_p = 0 \}$ denote the set of critical points of a 
smooth function $f:M \rightarrow \mathbb{R}$ on a  smooth $m$-dimensional manifold
$M$. A critical point $p \in Cr(f)$ is said to be \textbf{nondegenerate} if
and only if the Hessian $H_p(f)$ is nondegenerate.  The \textbf{index} $\lambda_p$
of a nondegenerate critical point $p$ is the dimension of the subspace of $T_pM$
where $H_p(f)$ is negative definite. If all the critical points of $f$ are
non-degenerate, then $f$ is called a \textbf{Morse function}.

If $f:M \rightarrow \mathbb{R}$ is a Morse function on a finite dimensional 
compact smooth Riemannian manifold $(M,g)$, then the \textbf{stable manifold} $W^s(p)$
and the \textbf{unstable manifold} $W^u(p)$ of a critical point $p \in Cr(f)$ are
defined to be
\begin{eqnarray*}
W^s(p) & = & \{ x\in M | \lim_{t \rightarrow \infty} \varphi_t(x) = p \}\\
W^u(p) & = & \{ x\in M | \lim_{t \rightarrow -\infty} \varphi_t(x) = p \}
\end{eqnarray*}
where $\varphi_t$ is the 1-parameter group of diffeomorphisms generated by
minus the gradient vector field, i.e. $-\nabla f$. The index of $p$ coincides 
with the dimension of $W^u(p)$.
The Stable/Unstable Manifold Theorem for a Morse Function says that the 
tangent space at $p$ splits as 
$$
T_pM = T^s_pM \oplus T_p^uM
$$
where the Hessian is positive definite on $T_p^sM \stackrel{def}{=} T_p W^s(p)$ and 
negative definite on $T_p^uM \stackrel{def}{=} T_p W^u(p)$.  Moreover, the
stable and unstable manifolds of $p$ are surjective images of smooth 
embeddings
\begin{eqnarray*}
E^s: T_p^sM & \rightarrow & W^s(p) \subseteq M\\
E^u: T_p^uM & \rightarrow & W^u(p) \subseteq M. 
\end{eqnarray*}
Hence, $W^s(p)$ is a smoothly embedded open disk of dimension $m - \lambda_p$, 
and $W^u(p)$ is a smoothly embedded open disk of dimension $\lambda_p$.

If the stable and unstable manifolds of a Morse function $f:M \rightarrow
\mathbb{R}$ all intersect transversally, then the function $f$ is called
\textbf{Morse-Smale}. Note that for a given Morse function $f:M \rightarrow
\mathbb{R}$ one can choose a Riemannian metric on $M$ so that $f$ is Morse-Smale
with respect to the chosen metric (see Theorem 2.20 of \cite{AbbLec}).
Moreover, if $f$ is Morse-Smale then $W(q,p) = W^u(q) \cap W^s(p)$ is an
embedded submanifold of $M$ of dimension $\lambda_q  - \lambda_p$, and when
$\lambda_q = \lambda_p + 1$ one can use Palis' $\lambda$-Lemma to prove 
that the number of gradient flow lines from $q$ to $p$ is finite.

If we assume that $M$ is oriented and we choose an orientation for each of the 
unstable manifolds of $f$, then there is an induced orientation on the stable 
manifolds. Thus, we can define an integer $n(q,p)$ associated to any 
two critical points $p$ and $q$ of relative index one by counting the 
number of gradient flow lines from $q$ to $p$ with signs determined by 
the orientations. If $M$ is not orientable, then we can still define $n(q,p)
\in \mathbb{Z}_2$ by counting the number of gradient flow lines from $q$ to $p$
$\text{mod }2$.

The \textbf{Morse-Smale-Witten chain complex} is defined to be the chain complex 
$(C_\ast(f),\partial_\ast)$ where $C_k(f)$ is the free abelian group
generated by the critical points $q$ of index $k$ (tensored with
$\mathbb{Z}_2$ when $n(q,p)$ is defined as an element of $\mathbb{Z}_2$) and the 
boundary operator $\partial_k:C_k(f) \rightarrow C_{k-1}(f)$ is given by
$$
\partial_k(q)\ \ = \sum_{p \in Cr_{k-1}(f)} n(q,p)p.
$$

\begin{theorem}[Morse Homology Theorem]\label{Morsehomology}
The pair $(C_\ast(f),\partial_\ast)$ is a chain complex.  If
$M$ is orientable and the boundary operator is defined
with $n(q,p) \in \mathbb{Z}$, then the homology
of $(C_\ast(f),\partial_\ast)$ is isomorphic to the singular homology 
$H_\ast(M;\mathbb{Z})$.  If the  boundary operator is defined with 
$n(q,p) \in \mathbb{Z}_2$, then the homology of $(C_\ast(f),\partial_\ast)$
is isomorphic to the singular homology $H_\ast(M;\mathbb{Z}_2)$.
\end{theorem}

\noindent
Note that the Morse Homology Theorem implies that the homology of
$(C_\ast(f),\partial_\ast)$ is independent of the Morse-Smale
function $f:M \rightarrow \mathbb{R}$, the Riemannian metric,
and the orientations.


\section{The Morse Inequalities}\label{MorseInequalities}

Let $M$ be a compact smooth manifold of dimension $m$. When $M$ is orientable
we define the $k^\text{th}$ \textbf{Betti number} of $M$, denoted $b_k$, to
be the rank of the $k^\text{th}$ homology group $H_k(M;\mathbb{Z})$ modulo its
torsion subgroup. When $M$ is not orientable we define $b_k = \text{dim }H_k(M;\mathbb{Z}_2)$.
Let $f:M \rightarrow \mathbb{R}$ be a Morse function on $M$, and let $\nu_k$ 
be the number of critical points of $f$ of index $k$ for all $k=0,\ldots ,m$.
As a consequence of the Morse Homology Theorem we have
$$
\nu_k \geq b_k
$$
for all $k=0,\ldots,m$ since $\nu_k = \text{rank } C_k(f)$ and
$H_k(M;\mathbb{Z})$ (or $H_k(M;\mathbb{Z}_2)$) is a quotient of $C_k(f)$.
These inequalities are known as the \textbf{weak Morse inequalities}.

\begin{definition}
The \textbf{Poincar\'e polynomial} \index{Poincar\'e polynomial} of $M$ 
is defined to be
$$
P_t(M) = \sum_{k=0}^m b_k t^k,
$$
and the \textbf{Morse polynomial} \index{Morse polynomial} of $f$ is defined 
to be
$$
M_t(f) = \sum_{k=0}^m \nu_k t^k.
$$
\end{definition}

\noindent
The Morse inequalities stated in the introduction are known as the
\textbf{strong Morse inequalities}. The strong Morse inequalities are
equivalent to the following \textbf{polynomial Morse inequalities}.
(For a proof see Lemma 3.43 of \cite{BanLec}.) It is this version of the Morse 
inequalities that we generalize to Morse-Bott functions.

\begin{theorem}[Polynomial Morse Inequalities]\label{Morsepolynomial}
For any Morse function $f:M \rightarrow \mathbb{R}$ on a smooth manifold $M$ 
we have
$$
M_t(f) = P_t(M) + (1+t) R(t)
$$
where $R(t)$ is a polynomial with non-negative integer coefficients.
That is, $R(t) = \sum_{k=0}^{m-1} r_k t^k$ where $r_k \in \mathbb{Z}$ 
satisfies $r_k \geq 0$ for all $k=0,\ldots ,m-1$.
\end{theorem}

\noindent
Although Theorem \ref{Morsepolynomial} and its proof are standard facts, we give here 
a detailed proof using the Morse-Smale-Witten chain complex because this proof gives an
explicit formula $r_k = \nu_{k+1} - z_{k+1}$ for the coefficients of the polynomial
$R(t)$ which we will use in the proof of Theorem \ref{MorseBottineq}.

\medskip\noindent
Proof of Theorem \ref{Morsepolynomial}:
Let $f:M \rightarrow \mathbb{R}$ be a Morse function on a finite dimensional compact 
smooth manifold $M$, and choose a Riemannian metric on $M$ for which $f$ is a Morse-Smale
function. Let $C_k(f)$ denote the $k^\text{th}$ chain group in the Morse-Smale-Witten chain
complex of $f$ (with coefficients in either $\mathbb{Z}$ or $\mathbb{Z}_2$), and let
$\partial_k:C_k(f) \rightarrow C_{k-1}(f)$ denote the $k^\text{th}$ Morse-Smale-Witten
boundary operator. The rank of $C_k(f)$ is equal to  the number $\nu_k$ of critical points of 
$f$ of index $k$, and by the Morse Homology Theorem (Theorem \ref{Morsehomology}) the rank of
$H_k(C_\ast(f),\partial_\ast) $ is equal to $b_k$, the $k^\text{th}$ Betti number of 
$M$ for all $k = 0,\ldots ,m$.

Let $z_k = \text{rank ker } \partial_k$ for all $k = 0,\ldots ,m$.
The exact sequence
$$
0 \rightarrow \text{ker }\partial_k \rightarrow C_k(f) 
\stackrel{\partial_k}{\rightarrow} \text{im }\partial_k
\rightarrow 0
$$
implies that $\nu_k = z_k + \text{rank im }\partial_k$ for all
$k = 0,\ldots, m$, and
$$
0 \rightarrow \text{im } \partial_{k+1} \rightarrow \text{ker }
\partial_k \rightarrow H_k(C_\ast(f),\partial_\ast)  \rightarrow 0
$$
implies that $b_k = z_k - \text{rank im } \partial_{k+1}$ for all
$k = 0,\ldots ,m$.  Hence,
\begin{eqnarray*}
M_t(f) - P_t(M) & = & \sum_{k=0}^m \nu_k t^k - \sum_{k=0}^m b_k t^k\\
& = & \sum_{k=0}^m (z_k + \text{rank im }\partial_k) t^k
    - \sum_{k=0}^m (z_k - \text{rank im }\partial_{k+1}) t^k\\
& = & \sum_{k=0}^m (\text{rank im }\partial_k + 
                    \text{rank im }\partial_{k+1}) t^k\\
& = & \sum_{k=0}^m (\nu_k - z_k + \nu_{k+1} - z_{k+1}) t^k\\
& = & \sum_{k=0}^m (\nu_k - z_k)t^k  + \sum_{k=0}^{m-1}(\nu_{k+1} - z_{k+1}) t^k\\
& = & t\sum_{k=1}^m (\nu_k - z_k)t^{k-1}  + \sum_{k=1}^{m}(\nu_k - z_k) t^{k-1} 
      \quad (\text{since } \nu_0 = z_0)\\
& = & (t+1)\sum_{k=1}^m (\nu_k - z_k) t^{k-1}.\\
\end{eqnarray*}
Therefore, $M_t(f) = P_t(M) + (1+t)R(t)$ where $R(t) = \sum_{k=0}^{m-1} (\nu_{k+1} - z_{k+1})t^k$.
Note that $\nu_{k+1} - z_{k+1} \geq 0$ for all $k=0,\ldots ,m-1$ because $z_{k+1}$ is the rank
of a subgroup of $C_{k+1}(f)$ and $\nu_{k+1}  =  \text{rank } 
C_{k+1}(f)$.

\proofend


\section{The Morse-Bott Inequalities}\label{MorseBottInequalities}

Let $f:M \rightarrow \mathbb{R}$ be a smooth function whose critical set
$\text{Cr}(f)$ contains a submanifold $C$ of positive dimension. 
Pick a Riemannian metric on $M$ and use it to split $T_\ast M|_C$
as
$$
T_\ast M|_C = T_\ast C \oplus \nu_\ast C
$$
where $T_\ast C$ is the tangent space of $C$ and $\nu_\ast C$ is the normal bundle of $C$.
Let $p \in C$, $V \in T_p C$, $W \in T_pM$, and let \index{Hessian} $H_p(f)$ 
be the Hessian of $f$ at $p$.  We have
$$
H_p(f)(V,W) = V_p \cdot (\tilde{W} \cdot f) = 0
$$
since $V_p \in T_pC$ and any extension of $W$ to a vector field $\tilde{W}$ satisfies 
$df(\tilde{W})|_C$ $= 0$.  Therefore, the Hessian $H_p(f)$ induces a
symmetric bilinear form $H_p^\nu(f)$ \index{$H_p^\nu(f)$ Hessian normal} on
$\nu_p C$.

\begin{definition}
A smooth function $f:M \rightarrow \mathbb{R}$ on a smooth manifold $M$ is 
called a \textbf{Morse-Bott function} if and only if the set of critical points
$\text{\rm Cr}(f)$ is a disjoint union of connected submanifolds and for each connected
submanifold $C \subseteq \text{\rm Cr}(f)$ the bilinear form $H_p^\nu(f)$ is non-degenerate 
for all $p \in C$.  
\end{definition}

\noindent
Often one says that the Hessian of a Morse-Bott function $f$ is 
non-degenerate in the direction normal to the critical submanifolds.

\medskip

For a proof of the following lemma see Section 3.5 of \cite{BanLec} or
\cite{BanApr}.

\begin{lemma}[Morse-Bott Lemma] \index{Morse-Bott Lemma}\label{MorseBottLemma}
Let $f:M \rightarrow \mathbb{R}$ be a Morse-Bott function and $C \subseteq
\text{\rm Cr}(f)$ a connected component.  For any $p \in C$ there is a local chart 
of $M$ around $p$ and a local splitting $\nu_\ast C = \nu_\ast^-C \oplus 
\nu_\ast^+C$, identifying a point $x\in M$ in its domain to $(u,v,w)$ where 
$u \in C$, $v \in \nu_\ast^-C$, $w \in \nu_\ast^+C$, such that within this chart
$f$ assumes the form
$$
f(x) = f(u,v,w) = f(C)- |v|^2 + |w|^2.
$$
\end{lemma}

\begin{definition}
Let $f:M \rightarrow \mathbb{R}$ be a Morse-Bott function on a finite 
dimensional smooth manifold $M$, and let $C$ be a critical
submanifold of $f$.  For any $p \in C$ let $\lambda_p$ denote the index of 
$H_p^\nu(f)$.  This integer is the dimension of $\nu_p^-C$ 
and is locally constant by the preceding lemma.  If $C$ is connected, then 
$\lambda_p$ is constant throughout $C$ and we call $\lambda_p = \lambda_C$ 
\index{$\lambda_C$ index of $C$} the \textbf{Morse-Bott index} \index{Morse-Bott index} 
of the connected critical submanifold $C$.
\end{definition}

\medskip

Let $f:M \rightarrow \mathbb{R}$ be a Morse-Bott function on a finite dimensional
compact smooth manifold, and assume that
$$
\text{Cr}(f) = \coprod_{j=1}^l C_j,
$$
where $C_1,\ldots ,C_l$ are disjoint connected critical submanifolds.  

\begin{definition}
The Morse-Bott polynomial of $f$ is defined to be
$$
MB_t(f) = \sum_{j=1}^l P_t(C_j) t^{\lambda_j}
$$
where $\lambda_j$ is the Morse-Bott index of the critical submanifold $C_j$
and $P_t(C_j)$ is the Poincar\'e polynomial of $C_j$.
\end{definition}

\noindent
Note: Clearly $MB_t(f)$ reduces to $M_t(f)$ when $f$ is a Morse function.  
Also, if $M$ or any of the critical submanifolds are not orientable, 
then $\mathbb{Z}_2$ coefficients are used to compute all the Betti numbers.

\smallskip
The following result is due to Bott.  In fact, Bott proved a more general version
of this result (without the assumption that the critical submanifolds are orientable)
using homology with local coefficients in an orientation bundle
\cite{BotNon} \cite{BotLec}.

\begin{theorem}[Morse-Bott Inequalities]\index{Morse-Bott Inequalities}\label{MorseBottineq}
Let $f:M \rightarrow \mathbb{R}$ be a Morse-Bott function on a finite
dimensional oriented compact smooth manifold, and assume that all the critical 
submanifolds of $f$ are orientable. Then there exists a polynomial 
$R(t)$ with non-negative integer coefficients such that
$$
MB_t(f) = P_t(M) + (1+t)R(t).
$$
If either $M$ or some of the critical submanifolds are not orientable, 
then this equation holds when $\mathbb{Z}_2$ coefficients are used to define
the Betti numbers.
\end{theorem}


\section{A perturbation technique and the proof of the main theorem}

To prove Theorem \ref{MorseBottineq} we will use the following perturbation 
technique based on \cite{AusMor}, the Morse-Bott Lemma, and a folk theorem 
proved in \cite{AbbLec}. This construction produces an explicit Morse-Smale 
function $h:M \rightarrow \mathbb{R}$ arbitrarily close to a given 
Morse-Bott function $f:M \rightarrow \mathbb{R}$ such that $h = f$ 
outside of a neighborhood of the critical set $\text{Cr}(f)$.

\smallskip
Let $T_j$ be a small tubular neighborhood around each connected
component $C_j\subseteq \text{Cr}(f)$ for all $j=1,\ldots ,l$
with local coordinates $(u,v,w)$ consistent with those from
the Morse-Bott Lemma (Lemma \ref{MorseBottLemma}). By ``small" we 
mean that each $T_j$ is contained in the union of the domains of the 
charts from the Morse-Bott Lemma, for $i\neq j$ we have $T_i \cap T_j = 
\emptyset$ and $f$ decreases by at least three times 
$\text{max}\{\text{var}(f,T_j)|\ j=1,\ldots, l\}$ along any gradient 
flow line from $T_i$ to $T_j$  where $\text{var}(f,T_j) = \sup\{f(x)|\ x 
\in T_j\} - \inf\{f(x)|\ x \in T_j\}$, and if $f(C_i) \neq f(C_j)$, then
$$
\text{var}(f,T_i) + \text{var}(f,T_j) < \left.\left.\frac{1}{3} \right|f(C_i) - f(C_j)\right|.
$$
Pick a Riemannian metric on $M$ such that the charts from the
Morse-Bott Lemma are isometries with respect to the standard
Euclidean metric on $\mathbb{R}^m$, and then pick positive Morse
functions $f_j:C_j \rightarrow \mathbb{R}$ that are Morse-Smale with
respect to the restriction of the Riemannian metric to $C_j$ for all
$j=1,\ldots ,l$. The Morse-Smale functions $f_j:C_j \rightarrow \mathbb{R}$
exist by the Kupka-Smale Theorem (see for instance Theorem 6.6. of \cite{BanLec}).

For every $j=1,\ldots ,l$ extend $f_j$ to a function on $T_j$ 
by making $f_j$ constant in the direction normal to $C_j$, i.e. $f_j$ is
constant in the $v$ and $w$ coordinates coming from the Morse-Bott 
Lemma. Let $\tilde{T}_j \subset T_j$ be a smaller tubular neighborhood
of $C_j$ with the same coordinates as $T_j$, and let $\rho_j$ 
be a smooth nonincreasing bump function which is constant in the $u$ coordinates,
equal to $1$ on $\tilde{T}_j$, and equal to $0$ outside of $T_j$. Choose $\varepsilon > 0$
small enough so that 
$$
\sup_{T_j - \tilde{T}_j} \varepsilon \|\nabla\rho_j f_j\| < \inf_{T_j - \tilde{T}_j} \|\nabla f\|
$$
for all $j=1,\ldots l$, and define
$$
h = f + \varepsilon \left( \sum_{j=1}^l \rho_j f_j \right).
$$
The function $h:M \rightarrow \mathbb{R}$ is a Morse function close to 
$f$, and the critical points of $h$ are exactly the critical points
of the $f_j$ for $j=1,\ldots ,l$. Moreover, if $p \in C_j$ is a critical point
of $f_j:C_j \rightarrow \mathbb{R}$ of index $\lambda_p^j$, then $p$ is a 
critical point of $h$ of index $\lambda_p^h = \lambda_j + \lambda_p^j$.

We now apply a well-known folk theorem (whose proof can be found in
Section 2.12 of \cite{AbbLec}) which says that we can perturb the Riemannian metric
on $M$ outside of the union of the tubular neighborhoods $T_j$ for $j=1,\ldots ,l$
so that $h$ satisfies the Morse-Smale transversality condition with respect to the
perturbed metric. In summary, we have achieved the following conditions.
\begin{enumerate}
\item The gradient $\nabla f = \nabla h$ outside of the union of the tubular 
      neighborhoods $T_j$ for $j=1,\ldots ,l$.  Moreover,
      the tubular neighborhood $T_j$ are chosen small enough so that the
      following condition holds: if $f(C_i) \leq f(C_j)$ for some $i\neq j$, 
      then there are no gradient flow lines of $h$ from $C_i$ to $C_j$. \label{fandh}
\item The charts from the Morse-Bott Lemma are isometries with respect
      to the metric on $M$ and the standard Euclidean metric on $\mathbb{R}^m$. \label{one}
\item In the local coordinates $(u,v,w)$ of a tubular neighborhood $T_j$
      we have $f = f(C)- |v|^2 + |w|^2$, $\rho_j$ depends only on the $v$ and $w$
      coordinates, and $f_j$ depends only on the $u$ coordinates. In particular, 
      $\nabla f \perp \nabla f_j$ on $T_j$ by the previous condition. \label{two}
\item The function $h = f + \varepsilon f_j$ on the tubular neighborhood $\tilde{T}_j$. \label{three}
\item The gradient $\nabla f$ dominates $\varepsilon \nabla \rho_j f_j$ on $T_j - \tilde{T}_j$. \label{four}
\item The functions $h:M \rightarrow \mathbb{R}$ and $f_j:C_j \rightarrow \mathbb{R}$
      satisfy the Morse-Smale transversality condition for all $j=1,\ldots ,l$.
\item For every $n=0,\ldots ,m$ we have the following description of the
$n^{\text{th}}$ Morse-Smale-Witten chain group of $h$ in terms of the Morse-Smale-Witten
chain groups of the $f_j$ for $j=1,\ldots ,l$. \label{decomposition}
$$
C_n(h) = \bigoplus_{\lambda_j + k = n} C_k(f_j)
$$
\end{enumerate}

\medskip
Now, let $\partial_\ast^h$ denote the Morse-Smale-Witten boundary operator of $h$, 
and let $\partial_\ast^{f_j}$ denote the Morse-Smale-Witten boundary 
operator of $f_j$ for $j=1,\ldots ,l$.  

\begin{lemma}\label{boundarylemma}
If $p,q \in C_j$ are critical points of $f_j:C_j \rightarrow \mathbb{R}$
of relative index one, then the coefficients $n(q,p)$ used to define
the Morse-Smale-Witten boundary operator are the same for $\partial_
\ast^h$ and $\partial_\ast^{f_j}$ (assuming in the case where $C_j$ 
is orientable that the appropriate orientation of $C_j$ has been chosen).
\end{lemma}

\smallskip\noindent
Proof: We have $h = f + \varepsilon f_j$ on $\tilde{T}_j$ by condition (\ref{three}).
This implies that $\nabla h = \nabla f + \varepsilon \nabla f_j$ on 
$\tilde{T}_j$, and $\nabla h = \varepsilon \nabla f_j$ on the critical submanifold
$C_j$. Moreover, by conditions (\ref{two}) and (\ref{four}) a gradient
flow line of $h$ cannot begin and end in $C_j$ unless the entire flow line 
is contained in $C_j$. Thus, if $p$ and $q$ are both in $C_j$ the flows
connecting them along the gradient flow lines of $h$ and of $f_j$ are the same, 
and the numbers $n(q,p)$ in the complex of $h$ and of $f_j$ are the 
same as long as the orientations are chosen appropriately.
\proofend

We now order the connected critical submanifolds $C_1,\ldots ,C_l$ by 
height, i.e. such that $f(C_i) \leq f(C_j)$ whenever $i \leq j$.
For $\alpha \in C_n(h)$ and $C_j$ a connected critical submanifold 
we denote by $\alpha_j$ the chain obtained from $\alpha$ by retaining 
only those critical points belonging to $C_j$. By condition (\ref{decomposition}), 
any $\alpha\neq 0$ can be written uniquely as:
$$
\alpha = \alpha_{j_1} + \cdots + \alpha_{j_r}
$$
where $j_1 < \cdots < j_r$ and $\alpha_{j_i}  \in C_{k_i}(f_{j_i}) - \{0\}$
where $\lambda_{j_i} + k_i = n$.  We will refer to $\alpha_{j_r}$ 
as the ``top part'' of  the chain $\alpha$.

\begin{corollary}\label{toppart}
If $\alpha$ is a non-zero element in $\text{ker } \partial_\ast^h$, then
the top part of $\alpha$ is a non-zero element in the kernel of $\partial_\ast^{f_j}$
for some $j$, i.e. using the above notation for $\alpha\in \text{ker } 
\partial_n^h$ we have $\alpha_{j_r} \in \text{ker }\partial_{k_r}^{f_{j_r}}$
where $\lambda_{j_r} + k_r = n$.
\end{corollary}

\smallskip\noindent
Proof: Let $\alpha = \sum_i n_iq_i \in \text{ker } \partial_n^h$, i.e.
$$
0 = \partial_n^h(\alpha)  =  \sum_i n_i \partial_n^h(q_i)
 =  \sum_i n_i \sum_{p \in Cr_{n-1}(h)} n(q_i,p)p
 =  \sum_{p \in Cr_{n-1}(h)} \left( \sum_i n_i n(q_i,p) \right) p.
$$
Then $\sum_i n_i n(q_i,p) = 0$ for all $p\in Cr_{n-1}(h)$, and in particular, for 
$p \in Cr_{n-1}(h) \cap C_{j_r}$. For any $p \in Cr_{n-1}(h) \cap C_{j_r}$
we have $f(q_i) \leq f(p)$ for all $q_i$ in the sum of $\alpha$ such that $q_i \not\in C_{j_r}$,
and thus there are no gradient flow lines of $h$ from $q_i$ to $p$ by condition
(\ref{fandh}).  Therefore, $n(q_i,p) = 0$ if $p \in  Cr_{n-1}(h) \cap
C_{j_r}$ and $q_i \not\in C_{j_r}$, i.e.
$$
0 = \sum_i n_i n(q_i,p) = \sum_{q_i \in C_{j_r}} n_i n(q_i,p)
$$
when $p \in  Cr_{n-1}(h) \cap C_{j_r}$. Lemma \ref{boundarylemma} 
then implies that $\partial_{k_r}^{f_{j_r}}(\alpha_{j_r}) = $
$$
\sum_{q_i \in C_{j_r}} n_i \partial_{k_r}^{f_{j_r}}(q_i)
= \sum_{q_i \in C_{j_r}} n_i \sum_{p \in Cr_{k_r-1}(f_{j_r})} n(q_i,p)p
= \sum_{p \in Cr_{k_r-1}(f_{j_r})} \left( \sum_{q_i \in C_{j_r}} n_i n(q_i,p) \right) p = 0,
$$
where $\lambda_{j_r} + k_r = n$ and $Cr_{k_r-1}(f_{j_r}) =  Cr_{n-1}(h) \cap C_{j_r}$. 
\proofend

\smallskip\noindent
Proof of Theorem \ref{MorseBottineq}: Let $f$, $h$, $f_j$ for $j=1,\ldots ,l$,
and the Riemannian metric on $M$ be as above. Let $M_t(f_j)$ denote 
the Morse polynomial of $f_j:C_j\rightarrow \mathbb{R}$, and let 
$c_j = \text{dim } C_j$ for all $j=1,\ldots ,l$. Note that the relation 
$\lambda_p^h = \lambda_j+ \lambda_p^j$ implies that
$$
M_t(h) = \sum_{j=1}^l M_t(f_j) t^{\lambda_j}.
$$
The polynomial Morse inequalities (Theorem \ref{Morsepolynomial}) say that
$$
M_t(h) = P_t(M) + (1+t)R_h(t)
$$
and
$$
M_t(f_j) = P_t(C_j) + (1+t)R_j(t)
$$
where $R_h(t)$ and $R_j(t)$ are polynomials with non-negative integer
coefficients for all $j=1,\ldots ,l$.
\begin{eqnarray*}
MB_t(f) & =  & \sum_{j=1}^l P_t(C_j) t^{\lambda_j}\\
& = & \left.\left.\sum_{j=1}^l \right( M_t(f_j) - (1+t) R_j(t) \right)t^{\lambda_j}\\
& = & \sum_{j=1}^l M_t(f_j)t^{\lambda_j} - (1+t) \sum_{j=1}^l R_j(t)t^{\lambda_j}\\
& = & M_t(h) - (1+t) \sum_{j=1}^l R_j(t)t^{\lambda_j}\\
& = & P_t(M) + (1+t) R_h(t) - (1+t) \sum_{j=1}^l R_j(t)t^{\lambda_j}
\end{eqnarray*}
In the proof of the polynomial Morse inequalities we saw that 
$$
R_j(t) = \sum_{k=1}^{c_j} (\nu_k^j - z_k^j)t^{k-1}
$$
where $\nu_k^j$ is the rank of the group $C_k(f_j)$ and $z_k^j$ is the rank
of the kernel of the boundary operator $\partial_k^{f_j}:C_k(f_j) \rightarrow C_{k-1}(f_j)$
in the Morse-Smale-Witten chain complex of $f_j:C_j\rightarrow \mathbb{R}$.
Hence,
$$
MB_t(f)= P_t(M) + (1+t)\sum_{n=1}^m (\nu_n^h - z_n^h) t^{n-1} - (1+t) 
\sum_{j=1}^l \sum_{k=1}^{c_j}(\nu_k^j - z_k^j) t^{\lambda_j+k-1}
$$
where $\nu_n^h$ denotes the rank of the chain group $C_n(h)$ and $z_n^h$ 
denotes the rank of the kernel of the Morse-Smale-Witten boundary operator
$\partial_n^h:C_n(h) \rightarrow C_{n-1}(h)$. Since the critical points of
$h$ coincide with the critical points of the functions $f_j$ for $j=1,\ldots l$,
and a critical point $p\in C_j$ of $f_j$ of index $\lambda_p^j$ is a critical
point of $h$ of index $\lambda_j + \lambda_p^j$ we have 
$$
\sum_{n=1}^m (\nu_n^h - z_n^h) t^{n-1} - \sum_{j=1}^l \sum_{k=1}^{c_j}(\nu_k^j - z_k^j) t^{\lambda_j+k-1}
= \sum_{j=1}^l \sum_{k=1}^{c_j}z_k^j t^{\lambda_j+k-1} - \sum_{n=1}^m z_n^h t^{n-1}.
$$
Therefore,
$$
MB_f(t) = P_t(M) + (1+t) R(t)
$$
where
$$
R(t) =  \sum_{j=1}^l \sum_{k=1}^{c_j}z_k^j t^{\lambda_j+k-1} - \sum_{n=1}^m z_n^h t^{n-1}
 = \sum_{n=1}^m \left( \sum_{\lambda_j + k = n} z_k^j - z_n^h \right) t^{n-1}.
$$

To see that $\sum_{\lambda_j + k = n} z_k^j \geq z_n^h$ for all $n=1,\ldots ,m$
we apply Corollary \ref{toppart} as follows.  If $z_n^h > 0$, then
we can choose a nonzero element $\beta_1 \in \text{ker } \partial_n^h$.
By Corollary \ref{toppart}, the ``top part'' of $\beta_1$ is a nonzero
element in $\text{ker }\partial_k^{f_j}$ for some $k$ and $j$ such that
$\lambda_j + k = n$. If $z_n^h = 1$, then we are done.  If $z_n^h > 1$,
then we can find an element $\beta_2 \in \text{ker } \partial_n^h$ that
is not in the group generated by $\beta_1$, and by adding a multiple of 
$\beta_1$ to $\beta_2$ if necessary, we can choose $\beta_2$ such that the
top part of $\beta_2$ is a nonzero element that is not in the subgroup 
generated by the top part of $\beta_1$. Continuing in this fashion we can 
find generators for $\text{ker }\partial^h_n$ whose top parts are
independent elements in 
$$
\bigoplus_{\lambda_j + k = n} \text{ker }\partial_k^{f_j},
$$
i.e.
$$
\sum_{\lambda_j + k = n} \left( \text{rank ker }\partial_k^{f_j} \right)
\geq \text{rank ker }\partial_n^h.
$$

\begin{flushright}
$\Box$
\end{flushright}

\smallskip\noindent
Acknowledgments:  We would like to thank William Minicozzi for several 
insightful conversations concerning Riemannian metrics and
gradient flows.


\def\cprime{$'$}

\end{document}